\documentclass[11pt,reqno]{amsart}
\usepackage[T1]{fontenc}
\usepackage{lmodern}
\usepackage{amssymb}
\usepackage[a4paper,margin=26mm]{geometry}
\usepackage[expansion=false]{microtype}
\usepackage{booktabs}
\usepackage{float}
\usepackage{caption}
\usepackage{xcolor}
\usepackage[hidelinks]{hyperref}
\hypersetup{
  pdftitle={The Extremal Number ex(Q7,C4) and Further Bounds for C4-Free Hypercube Subgraphs},
  pdfauthor={Zheli Zhou, Haitong Ma, Ruofeng Yang, Tianzhi Ren, Kuan Yang, Shuai Li},
  pdfsubject={Extremal graph theory on hypercubes},
  pdfkeywords={C4-free subgraphs, 4-cycle-free, hypercube, extremal number,
    edge Turan density, Erdos hypercube problem, flag algebras, SAT}
}

\numberwithin{equation}{section}
\newtheorem{theorem}{Theorem}[section]
\newtheorem{lemma}[theorem]{Lemma}
\newtheorem{proposition}[theorem]{Proposition}
\newcommand{\ex}{\operatorname{ex}}
\newcommand{\E}{\mathbb{E}}
\newcommand{\replayrevision}{c455c7788f3150160f33774d242f18cb91d93e7d}
\newcommand{\qeightrevision}{3d61e7ede69ebe898c89f1f8893a9afb861beae0}

\title[The extremal number \(\ex(Q_7,C_4)\)]
{The Extremal Number \(\ex(Q_7,C_4)\) and Further Bounds for
\(C_4\)-Free Hypercube Subgraphs}
\author[Z. Zhou]{Zheli Zhou}
\author[H. Ma]{Haitong Ma}
\author[R. Yang]{Ruofeng Yang}
\author[T. Ren]{Tianzhi Ren}
\author[K. Yang]{Kuan Yang}
\author[S. Li]{Shuai Li}
\address[Zheli Zhou, Haitong Ma, and Shuai Li]{Shanghai Innovation Institute, Shanghai, China}
\address[All authors]{Shanghai Jiao Tong University, Shanghai, China}
\email[Zheli Zhou]{zlz114514@sjtu.edu.cn}
\thanks{Shuai Li is the corresponding author.}

\begin{document}

\begin{abstract}
Let \(\ex(Q_n,C_4)\) be the maximum number of edges in a four-cycle-free
subgraph of the \(n\)-dimensional hypercube. We give computer-assisted
proofs that \(\ex(Q_7,C_4)=304\), that \(682\le\ex(Q_8,C_4)\le683\),
and that the edge Tur\'an density
\(\lim_{n\to\infty}\ex(Q_n,C_4)/(n2^{n-1})\) is less than
\(0.59088405\), improving Baber's \(0.60318\). In dimension seven, star and square counts reduce
the upper bound to local statements on six-dimensional faces,
verified by propositional refutations. In dimension eight, exact
sampling identities for faces, degrees, and rooted edge neighborhoods
yield a rational certificate excluding \(684\) or more edges.
For the asymptotic bound in Erd\H{o}s's hypercube problem, we use
flag algebra inequalities on coordinate faces and bound the
contribution of high-scoring seven-dimensional faces through their
rooted extension frequencies.
\end{abstract}

\maketitle

\section{Introduction}\label{sec:intro}

The hypercube \(Q_n\) has vertex set \(\{0,1\}^n\), with two vertices
adjacent when they differ in exactly one coordinate. Every four-cycle
alternates between two coordinate directions and hence is a coordinate
square. Following Erd\H{o}s~\cite[Problem 10]{erdos}, we consider the
extremal function
\[
  f(n)=\ex(Q_n,C_4)
       =\max\{e(G):G\subseteq Q_n,\ G\text{ is }C_4\text{-free}\}.
\]
Brass, Harborth, and Nienborg~\cite{bhn} gave bounds and constructions
for \(f(n)\). For \(2\le n\le4\), the exact values follow
from counting square--edge incidences and the constructions of
Derrida, Pomeau, Toulouse, and Vannimenus~\cite{dptv}.
The five- and six-dimensional values are given by Emamy-K., Guan, and
Rivera-Vega~\cite{egr} and Harborth and Nienborg~\cite{hn}, respectively.
Erd\H{o}s also asked whether density \(1/2+\varepsilon\) forces a
square for sufficiently large \(n\). Write
\[
  \rho(G)=\frac{e(G)}{n2^{n-1}}\quad(G\subseteq Q_n),
  \qquad
  \pi_\square(C_4)=\lim_{n\to\infty}\frac{f(n)}{n2^{n-1}}.
\]
The limit exists by face averaging (Lemma~\ref{lem:uniform-faces}).
Keeping the edges joining levels \(2j\) and \(2j+1\), where the level is
Hamming weight, gives a square-free graph with half the edges. Erd\H{o}s
conjectured that \(\pi_\square(C_4)=1/2\)~\cite{erdos}.
Thomason and Wagner~\cite{tw} obtained an upper bound of
approximately \(0.62256\). Balogh, Hu, Lidick\'y, and Liu~\cite{bhll}
obtained \(0.6068\), and Baber~\cite[Theorem 4.1]{baber} reached
\(0.60318\) using flag algebras on partially defined hypercubes.
Theorem~\ref{thm:asymptotic} improves this to \(0.59088405\).

An odd-square edge set meets every coordinate square in one or three
edges and is therefore \(C_4\)-free. In a fully frustrated signing,
every square has an odd number of negative edges, so the positive
edges form an odd-square set. The seven-dimensional construction of
Derrida et al.~\cite{dptv} gives \(304\) edges this way. For dimension
eight, Marinari, Parisi, and Ritort~\cite{mpr} reported a fully
frustrated ground state that corresponds to an odd-square set of
\(682\) edges~\cite[Section 6.3]{minamoto}. Minamoto~\cite[Theorem 1]{minamoto}
gives edge-list certificates and proves the odd-square optima in
dimensions six through eight. Minamoto also reports two SAT attempts and one
integer-programming attempt at the unrestricted seven-dimensional
problem, all inconclusive~\cite[Appendix A]{minamoto}. We resolve
Open Problem~1 in~\cite{minamoto}.

\begin{theorem}\label{thm:q7}
The maximum number of edges in a \(C_4\)-free subgraph of \(Q_7\) is
\(304\).
\end{theorem}

\begin{theorem}\label{thm:q8}
\(682\le f(8)\le 683\).
\end{theorem}

\begin{theorem}\label{thm:asymptotic}
\(\pi_\square(C_4)<0.59088405\).
\end{theorem}

Table~\ref{tab:extremal-values} lists the finite-dimensional values and bounds.

\begin{table}[htbp]
\caption{Extremal values and bounds through dimension eight.}
\label{tab:extremal-values}
\centering
\small
\begin{tabular}{@{}lrrrrrrrr@{}}
\toprule
Dimension \(n\) & 1 & 2 & 3 & 4 & 5 & 6 & 7 & 8\\
\(f(n)\) & 1 & 3 & 9 & 24 & 56 & 132 & 304 & \(682\)--\(683\)\\
\bottomrule
\end{tabular}
\end{table}

\subsection{Proof overview}

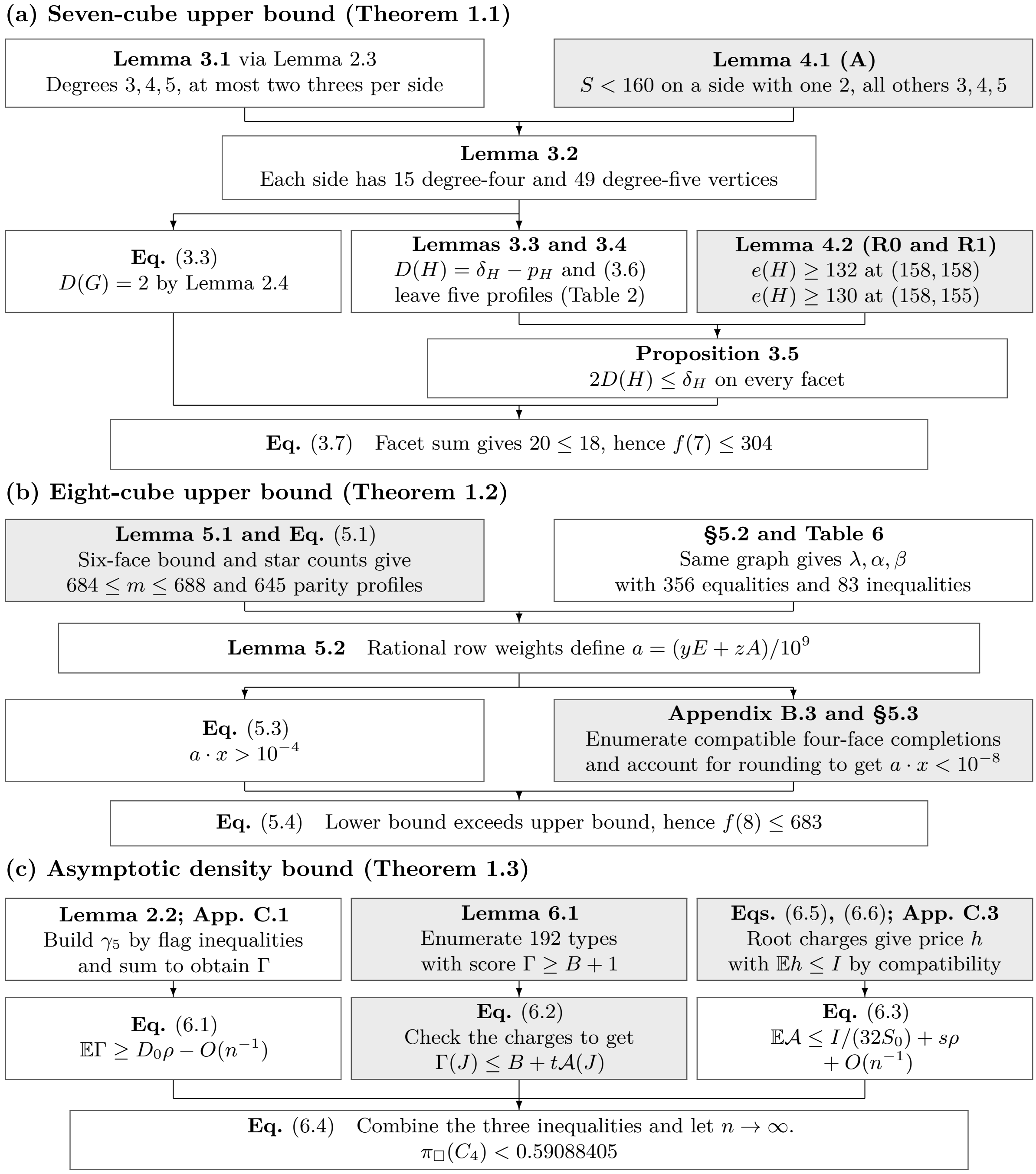
\begin{figure}[H]
\centering
\begingroup
\setlength{\unitlength}{1pt}
\setlength{\fboxsep}{3pt}
\setlength{\fboxrule}{0.45pt}
\newcommand{\proofnode}[6][white]{%
  \put(#2,#3){\makebox(0,0){\fcolorbox{black!65}{#1}{%
    \parbox[c][\dimexpr#5pt-2\fboxsep-2\fboxrule\relax][c]
      {\dimexpr#4pt-2\fboxsep-2\fboxrule\relax}
      {\centering\footnotesize #6}}}}}

\begin{picture}(446,205)
\put(0,197){\makebox(0,0)[l]{\small\bfseries
  (a) Seven-cube upper bound (Theorem~\ref{thm:q7})}}
\proofnode{104}{172}{208}{30}{%
  \textbf{Lemma~\ref{lem:initial}} via Lemma~\ref{lem:star}\\
  Degrees \(3,4,5\), at most two threes per side}
\proofnode[black!8]{342}{172}{208}{30}{%
  \textbf{Lemma~\ref{lem:root} (A)}\\
  \(S<160\) on a side with one \(2\), all others \(3,4,5\)}
\put(104,157){\line(0,-1){6}}
\put(342,157){\line(0,-1){6}}
\put(104,151){\line(1,0){238}}
\put(223,151){\vector(0,-1){6}}
\proofnode{223}{131}{258}{28}{%
  \textbf{Lemma~\ref{lem:degrees}}\\
  Each side has \(15\) degree-four and \(49\) degree-five vertices}
\put(223,117){\line(0,-1){6}}
\put(73,111){\line(1,0){150}}
\put(73,111){\vector(0,-1){7}}
\put(223,111){\vector(0,-1){7}}
\proofnode{73}{86}{146}{36}{%
  \textbf{Eq.~\eqref{eq:globalD}}\\
  \(D(G)=2\) by Lemma~\ref{lem:square}}
\proofnode{223}{86}{146}{36}{%
  \textbf{Lemmas~\ref{lem:local-identities} and~\ref{lem:opposite}}\\
  \(D(H)=\delta_H-p_H\) and \eqref{eq:opposite}\\
  leave five profiles (Table~\ref{tab:profiles})}
\proofnode[black!8]{373}{86}{146}{36}{%
  \textbf{Lemma~\ref{lem:rigid} (R0 and R1)}\\
  \(e(H)\ge132\) at \((158,158)\)\\
  \(e(H)\ge130\) at \((158,155)\)}
\put(223,68){\line(0,-1){5}}
\put(373,68){\line(0,-1){5}}
\put(223,63){\line(1,0){150}}
\put(309,63){\vector(0,-1){6}}
\proofnode{309}{44}{252}{26}{%
  \textbf{Proposition~\ref{prop:facet}}\\
  \(2D(H)\le\delta_H\) on every facet}
\put(73,68){\line(0,-1){40}}
\put(309,31){\line(0,-1){3}}
\put(73,28){\line(1,0){236}}
\put(223,28){\vector(0,-1){6}}
\proofnode{223}{11}{355}{22}{%
  \textbf{Eq.~\eqref{eq:q7-contradiction}}\quad
  Facet sum gives \(20\le18\), hence \(f(7)\le304\)}
\end{picture}\par\smallskip

\begin{picture}(446,160)
\put(0,153){\makebox(0,0)[l]{\small\bfseries
  (b) Eight-cube upper bound (Theorem~\ref{thm:q8})}}
\proofnode[black!8]{104}{124}{208}{36}{%
  \textbf{Lemma~\ref{lem:q8-gate} and Eq.~\eqref{eq:q8-budgets}}\\
  Six-face bound and star counts give\\
  \(684\le m\le688\) and \(645\) parity profiles}
\proofnode{342}{124}{208}{36}{%
  \textbf{\S\ref{sec:q8-sampling} and Table~\ref{tab:eight-rows}}\\
  Same graph gives \(\lambda,\alpha,\beta\)\\
  with \(356\) equalities and \(83\) inequalities}
\put(104,106){\line(0,-1){4}}
\put(342,106){\line(0,-1){4}}
\put(104,102){\line(1,0){238}}
\put(223,102){\vector(0,-1){5}}
\proofnode{223}{86}{400}{22}{%
  \textbf{Lemma~\ref{lem:q8-separation}}\quad
  Rational row weights define \(a=(yE+zA)/10^9\)}
\put(223,75){\line(0,-1){6}}
\put(104,69){\line(1,0){238}}
\put(104,69){\vector(0,-1){5}}
\put(342,69){\vector(0,-1){5}}
\proofnode{104}{46}{208}{36}{%
  \textbf{Eq.~\eqref{eq:q8-lower}}\\
  \(a\cdot x>10^{-4}\)}
\proofnode[black!8]{342}{46}{208}{36}{%
  \textbf{Appendix~\ref{app:eight-max} and \S\ref{sec:q8-separation}}\\
  Enumerate compatible four-face completions\\
  and account for rounding to get \(a\cdot x<10^{-8}\)}
\put(104,28){\line(0,-1){4}}
\put(342,28){\line(0,-1){4}}
\put(104,24){\line(1,0){238}}
\put(223,24){\vector(0,-1){4}}
\proofnode{223}{10}{355}{20}{%
  \textbf{Eq.~\eqref{eq:q8-margin}}\quad
  Lower bound exceeds upper bound, hence \(f(8)\le683\)}
\end{picture}\par\smallskip

\begin{picture}(446,137)
\put(0,130){\makebox(0,0)[l]{\small\bfseries
  (c) Asymptotic density bound (Theorem~\ref{thm:asymptotic})}}
\proofnode{73}{100}{146}{36}{%
  \textbf{Lemma~\ref{lem:uniform-faces}; App.~\ref{app:five-score}}\\
  Build \(\gamma_5\) by flag inequalities\\
  and sum to obtain \(\Gamma\)}
\proofnode[black!8]{223}{100}{146}{36}{%
  \textbf{Lemma~\ref{lem:high-score-catalogue}}\\
  Enumerate \(192\) types\\
  with score \(\Gamma\ge B+1\)}
\proofnode[black!8]{373}{100}{146}{36}{%
  \textbf{Eqs.~\eqref{eq:asymp-root}, \eqref{eq:asymp-price}; App.~\ref{app:charge}}\\
  Root charges give price \(h\)\\
  with \(\E h\le I\) by compatibility}
\put(73,82){\vector(0,-1){7}}
\put(223,82){\vector(0,-1){7}}
\put(373,82){\vector(0,-1){7}}
\proofnode{73}{57}{146}{36}{%
  \textbf{Eq.~\eqref{eq:asymp-score}}\\
  \(\E\Gamma\ge D_0\rho-O(n^{-1})\)}
\proofnode[black!8]{223}{57}{146}{36}{%
  \textbf{Eq.~\eqref{eq:asymp-point}}\\
  Check the charges to get\\
  \(\Gamma(J)\le B+t\mathcal A(J)\)}
\proofnode{373}{57}{146}{36}{%
  \textbf{Eq.~\eqref{eq:asymp-charge}}\\
  \(\E\mathcal A\le I/(32S_0)+s\rho\)\\
  \({}+O(n^{-1})\)}
\put(73,39){\line(0,-1){8}}
\put(373,39){\line(0,-1){8}}
\put(73,31){\line(1,0){300}}
\put(223,39){\vector(0,-1){13}}
\proofnode{223}{13}{390}{26}{%
  \textbf{Eq.~\eqref{eq:asymp-result}}\quad
  Combine the three inequalities and let \(n\to\infty\).\\
  \(\pi_\square(C_4)<0.59088405\)}
\end{picture}
\endgroup
\caption{Proof dependencies. Gray boxes use finite computation. The finite
lower bounds follow from the constructions discussed in Section~\ref{sec:intro}.}
\label{fig:proof-routes}
\end{figure}

Figure~\ref{fig:proof-routes} traces the dependencies of the three
upper bounds. Each arrow points from an input to the step that uses it.

\begin{samepage}
\emph{The seven-cube.}
A hypothetical \(305\)-edge example \(G\) has degree sum \(305\)
on each parity side.
The linear bound \(\binom d3\ge6d-20\) gives at least \(550\) stars
per side, while each side can support at most \(560\). The resulting
slack of ten gives initial degree restrictions. Summing refined
star counts over the six-dimensional facets and invoking
Obstruction A, which excludes a side with one vertex of degree two and
a saturated star count (Lemma~\ref{lem:root}), then forces all degrees
to lie in \(\{4,5\}\).
\par
\end{samepage}

For a graph \(X\), let \(D(X)\) be its number of empty squares minus
its number of squares containing exactly two opposite edges. The
degree sequence determines \(D(G)=2\). For a six-dimensional
restriction \(H\), write \(\delta_H=132-e(H)\). Comparing \(H\) with
its opposite restriction leaves five candidate violations of
\(2D(H)\le\delta_H\). The star-count bounds R0 and R1 exclude
these violations when the two parity star counts are \((158,158)\)
or \((158,155)\), respectively (Lemma~\ref{lem:rigid}).

Each square lies in five facets and each edge in six, so summing over
the fourteen facet restrictions gives
\[
  20=10D(G)=2\sum_H D(H)\le\sum_H\delta_H=18,\qquad\text{a contradiction.}
\]

\emph{The eight-cube.}
An inclusion-maximal counterexample has between \(684\) and \(688\)
edges by a one-sided \(Q_6\) inequality and face counting. Its
four-face distribution, two parity-side degree profiles, and ambient
edge-neighborhood distribution obey exact sampling identities.
A rational combination of these identities and local inequalities
gives a lower bound on an averaged score. Exact maximization over the
three finite domains gives a strictly smaller upper bound.

\emph{The density bound.}
The first step is to add corrections whose expectations are
asymptotically nonnegative to the edge density of a random face.
The resulting seven-face score can be large on a finite list of
types. We bound their effect by a nonnegative charge on rooted
extensions and control its expectation using smaller faces that
share the same edges.

Section~\ref{sec:prelims} gives the face-counting identities.
Sections~\ref{sec:q7} and~\ref{sec:computation} prove the
result for \(Q_7\). Sections~\ref{sec:q8} and~\ref{sec:asymptotic}
give the eight-dimensional and asymptotic bounds. The appendices
describe the finite domains and the verification of the certificates.

\section{Coordinate faces and counting identities}\label{sec:prelims}

We treat subgraphs as spanning by allowing isolated vertices. Let
\(P_0\) and \(P_1\) be the even and odd parity sides, determined by the
parity of the Hamming weight. Coordinate indices run from \(0\) to
\(n-1\). A \(k\)-dimensional coordinate face is obtained
by choosing \(k\) free coordinates and fixing the remaining coordinates.
We use \(G[F]\) for the restriction of a graph \(G\) to a face \(F\).
An \((n-1)\)-dimensional coordinate face of \(Q_n\) is a \emph{facet}.
In the reduction, \(F\) denotes an ambient facet and \(H=G[F]\) its
restriction. The notation \(\sum_H\) ranges over all facet restrictions,
indexed by their ambient facets. Write \(d_G(v)\) for degree and \(e(G)\)
for edge count. For a statement \(E\), the Iverson bracket \([E]\) is
\(1\) if \(E\) holds and \(0\) otherwise.

\begin{lemma}[Face multiplicities]\label{lem:multiplicities}
There are \(\binom nk 2^{n-k}\) coordinate \(k\)-faces of \(Q_n\), and
each coordinate \(r\)-face lies in \(\binom{n-r}{k-r}\) of them for
\(r\le k\le n\).
\end{lemma}

\begin{proof}
Choose the \(k\) free coordinates and assign values to the remaining
\(n-k\). To extend an \(r\)-face to a \(k\)-face, choose \(k-r\) of
its fixed coordinates to be freed.
\end{proof}

\begin{lemma}[Uniform face sampling]\label{lem:uniform-faces}
If \(F\) is a uniform \(k\)-face of \(Q_n\), with \(1\le k\le n\),
then \(\E\rho(G[F])=\rho(G)\).
If a uniform \(r\)-face is sampled inside a uniform \(k\)-face, the
result is a uniform \(r\)-face of \(Q_n\).
\end{lemma}

\begin{proof}
The number of \(k\)-faces containing each fixed edge is constant by
Lemma~\ref{lem:multiplicities}, and all \(k\)-faces have the same
number of ambient edges. This gives the density identity by
double counting. The same constant-multiplicity argument applies
to a fixed \(r\)-face. In particular, \(f(n)/(n2^{n-1})\) is
nonincreasing, since every restriction of a \(C_4\)-free graph
is \(C_4\)-free.
\end{proof}

A \emph{full three-star} consists of three retained edges incident
with a common vertex, its root. The three directions determine a
unique coordinate \(3\)-face containing the star. For a parity side
\(P\), write
\[
  S_P(G)=\sum_{v\in P}\binom{d_G(v)}{3}.
\]

\begin{lemma}[Star capacity]\label{lem:star}
For a \(C_4\)-free graph \(G\subseteq Q_n\), each coordinate \(3\)-face contains at most
one full three-star rooted on a specified parity side. Consequently
\[
  S_P(G)\le \binom{n}{3}2^{n-3}.
\]
The difference counts the \(3\)-faces \emph{missing a \(P\)-star}.
\end{lemma}

\begin{proof}
Two distinct same-parity vertices of a \(3\)-face differ in exactly
two coordinates. If both were roots of full stars in the face, the four
edges in those two directions incident with the two roots would form a
\(C_4\). Summing over the \(3\)-faces proves the bound and the counting
interpretation.
\end{proof}

Write \(a_0(G)\) for the number of empty coordinate squares and
\(a_{\mathrm{opp}}(G)\) for the number with exactly two retained edges
that are opposite. Define \(D(G)=a_0(G)-a_{\mathrm{opp}}(G)\).

\begin{lemma}[Square identity]\label{lem:square}
If \(G\subseteq Q_n\) is \(C_4\)-free, then
\begin{equation}\label{eq:square}
  \sum_{v}d_G(v)^2
    =2n\,e(G)-2\binom n2 2^{n-2}+2D(G).
\end{equation}
\end{lemma}

\begin{proof}
Every adjacent pair of retained edges belongs to a unique square. In a
square with \(k\) retained edges, the number of adjacent pairs is
\(k-1+[k=0]-[\text{two opposite edges}]\),
as is seen by considering \(k=0,1,2,3\). Each retained edge lies in
\(n-1\) squares. Summing over squares gives
\[
  \sum_v\binom{d_G(v)}{2}=(n-1)e(G)-\binom{n}{2}2^{n-2}+D(G).
\]
Now use \(d^2=2\binom d2+d\) and \(\sum_v d_G(v)=2e(G)\).
\end{proof}

\section{The seven-cube: reduction to local obstructions}\label{sec:q7}

Deleting edges preserves \(C_4\)-freeness, so for the upper bound it
suffices to rule out the assumption
\begin{equation*}
  G\subseteq Q_7\text{ is \(C_4\)-free and }e(G)=305.
  \tag{\(\star\)}\label{hyp:counterexample}
\end{equation*}
Throughout this section, \(G\) satisfies \eqref{hyp:counterexample}.
The two finite lemmas used in the reduction are stated and proved in
Section~\ref{sec:computation}.

\subsection{Global degree restrictions}

Throughout this subsection, counts \(a,b,c\) refer to one chosen
parity side \(P\).

\begin{lemma}[Initial degree restrictions]\label{lem:initial}
On each parity side of \(G\), the degrees belong to \(\{3,4,5\}\)
and at most two vertices have degree three.
\end{lemma}

\begin{proof}
\emph{Degree penalty.}
Let \(\omega(d)=\binom d3-6d+20\). For \(d=0,\ldots,7\) the values are
\[
  20,\quad 14,\quad 8,\quad 3,\quad 0,\quad 0,\quad 4,\quad 13 .
\]
Since each side has \(64\)
vertices with degree sum \(305\), Lemma~\ref{lem:star} gives
\[
  \sum_{v\in P}\omega(d_G(v))=S_P(G)-550\le10 .
\]
In particular, degrees \(0,1,7\) are excluded. Let \(a,b,c\) count
vertices of degrees \(2,3,6\) in \(P\), respectively, and let \(n_4(P)\)
count vertices of degree four in \(P\). The penalty
bound and the degree sum yield
\begin{equation}\label{eq:penalty}
  8a+3b+4c\le10,\qquad n_4(P)=15-3a-2b+c.
\end{equation}

\emph{Facet bound.}
For each restriction \(H=G[F]\), put \(\delta_H=132-e(H)\) and
let \(A_H,B_H,C_H\) count vertices of local degree \(1,2,6\), respectively,
in \(P\cap V(H)\). Write \(g(d)=(d-4)(d-5)/2\). For the possible
local degrees \(1\le d\le6\), one has
\[
  \binom d3=3g(d)+6d-20
      -4[d=1]-[d=2]+[d=6] .
\]
The local degree sum on one side is \(e(H)\), and its star count is
at most \(160\). Rearranging and rounding down gives
\begin{equation}\label{eq:floor}
  \sum_{v\in P\cap V(H)}g(d_H(v))
  \le 2+2\delta_H+
      \left\lfloor\frac{2+4A_H+B_H-C_H}{3}\right\rfloor .
\end{equation}

\emph{Summation over facets.}
A vertex of degree \(d\) occurs in seven facets, with local degree
\(d-1\) in \(d\) of them and \(d\) in the remaining \(7-d\).
Its contribution to the sum of the left-hand sides is
\(d\,g(d-1)+(7-d)g(d)\), taking values \(27,13,4,0,1\)
for \(d=2,3,4,5,6\). By \eqref{eq:penalty} the total is
\(60+15a+5b+5c\). Each edge belongs to six facets, so
\(\sum_H\delta_H=14\cdot132-6\cdot305=18\).

If \(a=0\), then \(A_H=0\). Since \(C_H\ge0\) and \(B_H\) is a
nonnegative integer,
\[
  \left\lfloor\frac{2+B_H-C_H}{3}\right\rfloor
  \le \left\lfloor\frac{2+B_H}{3}\right\rfloor
  \le B_H .
\]
Each global degree-three vertex contributes to
\(B_H\) in three facets, giving \(\sum_H B_H=3b\).
If \(a=1\), then \eqref{eq:penalty} gives \(b=c=0\). The sole
degree-two vertex contributes \(A_H=1\) in two facets and
\(B_H=1\) in five others, and the sum of the floors is \(9\).
These cases exhaust \(a\) by \eqref{eq:penalty}. Summing
\eqref{eq:floor} gives
\[
  60+15a+5b+5c\le 64+9a+3b,
  \qquad 6a+2b+5c\le4 .
\]
Thus \(a=c=0\) and \(b\le2\) on the chosen side. The argument
applies to either parity side.
\end{proof}

\begin{lemma}[Degree sequence]\label{lem:degrees}
Each parity side of \(G\) has \(15\) vertices of degree four and
\(49\) of degree five.
\end{lemma}

\begin{proof}
Fix one parity side \(P\). By Lemma~\ref{lem:initial} the number \(b\)
of degree-three vertices on this side belongs to \(\{0,1,2\}\).

If \(b=2\), then \(A_H=C_H=0\) everywhere. The bounds in the proof of
Lemma~\ref{lem:initial} give
\[
  70=60+5b
  =\sum_H\sum_{v\in P\cap V(H)}g(d_H(v))
  \le 64+\sum_H\left\lfloor\frac{2+B_H}{3}\right\rfloor
  \le 64+\sum_H B_H=70 .
\]
The two inequalities sum nonnegative slacks, so both are tight in
each facet. In particular, \(\lfloor(2+B_H)/3\rfloor=B_H\) for every
facet, which forces \(B_H\in\{0,1\}\).
Since \(\sum_H B_H=3b=6\), six facets have \(B_H=1\). In each such
facet equality in \eqref{eq:floor} implies a star count of \(160\).
Indeed, \(2+B_H=3\) is divisible by three, so the rounding step
loses nothing there. Since \(B_H=1\), exactly one vertex on this side
drops from degree three to degree two. All other local degrees lie in
\(\{3,4,5\}\), since the global degrees are in \(\{3,4,5\}\) and
each can drop by at most one. This contradicts
Lemma~\ref{lem:root}.

If \(b=1\), the unique degree-three vertex of \(P\) belongs to seven
facets. It drops to degree two in three and retains degree three in
four. The seven facets omitting the vertex and the four in which it
retains its degree all have their local degrees on \(P\) in
\(\{3,4,5\}\). Each of these degrees contributes \(1\) modulo \(3\)
to the star count. The eleven corresponding counts are congruent to
\mbox{\(32\equiv2\pmod{3}\)}, so the capacity of \(160\) gives an upper
bound of \(158\). In the other three facets the degree-two vertex
contributes zero, and the counts are congruent to
\mbox{\(31\equiv1\pmod{3}\)}. Lemma~\ref{lem:root} excludes \(160\),
leaving an upper bound of \(157\).
But \(S_P(G)=6\cdot305-20\cdot64+3=553\).
Each full three-star occurs in four facets by
Lemma~\ref{lem:multiplicities}. The sum of the facet star counts is
then \(4\cdot553=2212\), exceeding \(3\cdot157+11\cdot158=2209\).

Consequently \(b=0\) on either side. The degree sum \(305\) on a
side of size \(64\) gives the stated counts.
\end{proof}

The degree sequence in Lemma~\ref{lem:degrees} gives
\(\sum_v d_G(v)^2=2(15\cdot4^2+49\cdot5^2)=2930\).
Using \eqref{eq:square} with \(n=7\) and \(e(G)=305\) yields
\begin{equation}\label{eq:globalD}
  D(G)=2.
\end{equation}

\subsection{Opposite facets and the facet inequality}

In each restriction \(H=G[F]\), every local degree is in
\(\{3,4,5\}\). For \(i=0,1\) write \(b_{H,i}\) and \(S_{H,i}\)
for the degree-three count and star count on its side inherited from
\(P_i\). Set
\[
  \delta_H=132-e(H),\qquad p_{H,i}=(158-S_{H,i})/3,\qquad
  p_H=p_{H,0}+p_{H,1}.
\]
For \(d=3,4,5\), the identity
\[
  \binom{d}{3}=6d-20+3[d=3]
\]
follows by evaluating both sides at those degrees, giving \(1,4,10\),
respectively. Summing over the \(32\) vertices gives
\begin{equation}\label{eq:local-stars}
  S_{H,i}=6e(H)-640+3b_{H,i}.
\end{equation}
Equation~\eqref{eq:local-stars} gives \(S_{H,i}\equiv2\pmod{3}\).
The capacity bound \(S_{H,i}\le160\) from Lemma~\ref{lem:star} therefore
implies \(S_{H,i}\le158\), so \(p_{H,i}\) is a nonnegative integer.

\begin{lemma}[Local identities]\label{lem:local-identities}
Every facet restriction \(H\) satisfies
\begin{equation}\label{eq:local-identities}
  \delta_H\ge0,\qquad
  b_{H,i}=2+2\delta_H-p_{H,i},\qquad D(H)=\delta_H-p_H .
\end{equation}
\end{lemma}

\begin{proof}
Substituting the definitions of \(p_{H,i}\) and \(\delta_H\) into
\eqref{eq:local-stars} gives the identity for \(b_{H,i}\).
For degrees \(3,4,5\), the polynomial \((d-4)(d-5)\) equals
\(2[d=3]\). Summing and using \eqref{eq:square} in
dimension six gives
\[
  b_{H,0}+b_{H,1}=400-3e(H)+D(H)=4+3\delta_H+D(H).
\]
Comparing with \(b_{H,0}+b_{H,1}=4+4\delta_H-p_H\) yields the
identity for \(D(H)\).

To prove nonnegativity of \(\delta_H\), note that
\(6e(H)-640\le S_{H,i}\le158\), so \(e(H)\le133\).
If \(e(H)=133\), both parity star counts are \(158\), and
\eqref{eq:local-stars} forces \(b_{H,0}=b_{H,1}=0\).
All degrees in \(H\) are then four or five. On restriction to any
five-dimensional face \(K\) of the ambient six-cube, the local degrees
belong to \(\{3,4,5\}\), since fixing one more coordinate removes at
most one incident edge.
Each of the \(16\) vertices on the chosen side of \(K\) therefore
contributes \(1\) modulo \(3\) to its star count \(S\).
There are \(40\) coordinate \(3\)-faces in \(K\), so the number of its
faces missing a star on that side is \(40-S\equiv0\pmod{3}\).
Star status in a coordinate \(3\)-face depends only on its own edges.
The missing faces in \(K\) are therefore among the two missing faces
on that side of \(H\). It follows that every \(K\) has zero missing
faces, which is impossible since each missing \(3\)-face of \(H\)
lies in three \(5\)-facets. This proves \(e(H)\le132\) for the
restrictions in the reduction without using the known value
\(f(6)=132\)~\cite{hn}.
\end{proof}

\begin{lemma}[Opposite facets]\label{lem:opposite}
If \(H\) and \(H'\) are restrictions to opposite facets, then
\begin{equation}\label{eq:opposite}
  3\delta_H+\delta_{H'}\le11+p_H .
\end{equation}
\end{lemma}

\begin{proof}
Let \(L_{H,i}\) count globally degree-four vertices of \(P_i\)
in \(H\), and let \(L_H=L_{H,0}+L_{H,1}\). Let \(t\) count
retained edges between \(P_0\cap V(H)\) and \(P_1\cap V(H')\).
The global degree sums on these two sets are
\[
  e(H)+t=160-L_{H,0},
  \qquad e(H')+t=160-L_{H',1}.
\]
Subtract and use \(L_{H,1}+L_{H',1}=15\) to obtain
\(\delta_H-\delta_{H'}=L_H-15\).
Each vertex of local degree three has global degree four, so
\(L_H\ge b_{H,0}+b_{H,1}=4+4\delta_H-p_H\).
Rearranging gives \eqref{eq:opposite}.
\end{proof}

\begin{proposition}[Facet inequality]\label{prop:facet}
Every facet restriction \(H\) of \(G\) satisfies \(2D(H)\le\delta_H\).
\end{proposition}

\begin{proof}
If \(2D(H)>\delta_H\), then \(2p_H<\delta_H\) by
\eqref{eq:local-identities}. The nonnegativity of \(\delta_{H'}\)
and \eqref{eq:opposite} give
\[
  3\delta_H\le11+p_H\le11+\frac{\delta_H-1}{2},
  \qquad 5\delta_H\le21 .
\]
Thus \(1\le\delta_H\le4\). For each such \(\delta_H\), the integer \(p_H\)
satisfies
\[
  \max\{0,3\delta_H-11\}\le p_H
  \le \left\lfloor\frac{\delta_H-1}{2}\right\rfloor .
\]
The lower bound excludes \(p_H=0\) when \(\delta_H=4\). This leaves the
five rows of Table~\ref{tab:profiles}, up to exchanging the parity
sides. The rows with \(p_H=0\) have star counts
\((158,158)\), and the rows with \(p_H=1\) have star counts
\((158,155)\) up to order. Lemma~\ref{lem:rigid} rules out every row.
\end{proof}

\begin{table}[H]
\centering
\caption{Five candidate violations of the facet inequality.}\label{tab:profiles}
\small
\begin{tabular}{@{}rrrccc@{}}
\toprule
\(\delta_H\) & \(p_H\) & \(e(H)\) &
\((b_{H,0},b_{H,1})\) & \((S_{H,0},S_{H,1})\) & Obstruction\\
\midrule
1 & 0 & 131 & \((4,4)\) & \((158,158)\) & \(\mathrm{R0}\)\\
2 & 0 & 130 & \((6,6)\) & \((158,158)\) & \(\mathrm{R0}\)\\
3 & 0 & 129 & \((8,8)\) & \((158,158)\) & \(\mathrm{R0}\)\\
3 & 1 & 129 & \((7,8)\) & \((155,158)\) & \(\mathrm{R1}\)\\
4 & 1 & 128 & \((9,10)\) & \((155,158)\) & \(\mathrm{R1}\)\\
\bottomrule
\end{tabular}
\end{table}

\subsection{Completion of the proof}

\begin{proof}[Proof of Theorem~\ref{thm:q7}]
Each square lies in five \(6\)-facets and makes the same contribution
to \(D\) in each restriction. Each edge lies in six \(6\)-facets.
Summing Proposition~\ref{prop:facet} over all fourteen facets and using
\eqref{eq:globalD} gives
\begin{equation}\label{eq:q7-contradiction}
  20=2\cdot5D(G)=\sum_H 2D(H)
    \le \sum_H\delta_H
    =14\cdot132-6\cdot305=18 ,
\end{equation}
contradicting \eqref{hyp:counterexample}. Minamoto's
\(304\)-edge witness~\cite[Proposition 4]{minamoto}\footnote{\raggedright Edge list at
\url{https://github.com/minamominamoto/c4free-hypercube/blob/d55a01b9b7a64e3b81baca72c216590a9b1b1c41/selected_edges_best.json}.}
gives the matching lower bound. It has squares with two retained edges.
\end{proof}

\section{Local obstructions for the seven-cube}\label{sec:computation}

\subsection{Finite statements}

\begin{lemma}[Obstruction A]\label{lem:root}
There is no \(C_4\)-free subgraph of \(Q_6\) with a parity side that
has one vertex of degree two, has degree in \(\{3,4,5\}\) at each of
its other \(31\) vertices, and has star count \(160\). The opposite
parity side has no degree constraint.
\end{lemma}

\begin{proof}
Translate the degree-two vertex to the origin and permute coordinates
so that its two retained directions are \(0\) and \(1\). The
specified side is then \(P_0\). Use a Boolean
variable for each of the \(192\) edges and an indicator for each of
the \(640\) candidate full stars rooted on the even side. For each
square add a clause forbidding its four edges from all being retained.
Set each star indicator equal to the conjunction of its three edges.
The star-count hypothesis and Lemma~\ref{lem:star} require exactly
one even-rooted star in each of the \(160\) three-dimensional faces.
Fix the six incident edge variables at the origin and impose the
stated degree bounds at the remaining \(31\) even vertices.

A graph in the statement would satisfy this formula, with the star
indicators set by its edges. The resulting CNF has \(832\) variables
and \(4{,}422\) clauses. Section~\ref{sec:replay} records a verified
refutation.
\end{proof}

\begin{lemma}[Two star-count bounds]\label{lem:rigid}
Let \(H\subseteq Q_6\) be \(C_4\)-free, with all degrees in
\(\{3,4,5\}\).
\begin{enumerate}
\item[\(\mathrm{R0}\).] If the two parity star counts are \(158,158\),
  then \(e(H)\ge132\).
\item[\(\mathrm{R1}\).] If the two parity star counts are \(158,155\)
  in either order, then \(e(H)\ge130\).
\end{enumerate}
\end{lemma}

\begin{proof}
By \eqref{eq:local-stars}, on a side with \(S=158\) and \(b\)
degree-three vertices one has \(b=266-2e(H)\).
A counterexample to \(\mathrm{R0}\) has \(b\ge4\)
on that side, and a counterexample to \(\mathrm{R1}\) has \(b\ge8\)
on the side with \(S=158\).

\(Q_6\) has \(160\) coordinate \(3\)-faces, so this side misses
a star in exactly two of them. Normalize the first such face to the
face with free coordinates \(0,1,2\) and all other coordinates zero.
To make the side of interest even, use a vertex from that side in the
first missing face as the translation vector.
The stabilizer of this face preserving even parity has ten orbits on
the remaining faces, derived in Appendix~\ref{app:orbits}. For each
orbit, require the stated two missing even stars, the appropriate
number of missing odd stars (two or five), the degree bounds, and the
lower bound on the even-side \(b\). In R1 the side with \(S=155\)
has its five missing stars counted without fixing their positions or
imposing a separate degree-three count. Each counterexample therefore yields a satisfying
assignment of one of the ten CNF inputs for the corresponding assertion.
Appendix~\ref{app:encoding} gives the clauses, and
Section~\ref{sec:replay} records the checks of all twenty refutations.
\end{proof}

The bound in R0 is attained. Minamoto's
\(132\)-edge odd-square witness~\cite[Section 6.3]{minamoto}\footnote{\raggedright Edge list at
\url{https://github.com/minamominamoto/c4free-hypercube/blob/918025df85b509c4311675e0d1ba65f26f23d98d/q6_odd_square_132.json}.}
has two degree-three, twenty-four degree-four, and six degree-five
vertices on each side. Thus its two star counts are \(158\).

\subsection*{Seven-dimensional verification data}

\url{https://github.com/zhelizhou/hypercube-c4} is the companion repository.
The frozen \path{q7/} inputs and original replay records are at commit
\href{https://github.com/zhelizhou/hypercube-c4/tree/13c4ac5b2269383d36e5fb0da303f0a5cd2d1b03/q7}
{\texttt{13c4ac5b2269}}. The CNF generator reproduces their exact bytes
using the reconstruction routine in the encoding audit.
The release containing all \(21\) traces is at
\par\nopagebreak\noindent
\url{https://github.com/zhelizhou/hypercube-c4/releases/tag/q7-artifacts-v0.1}.

\subsection{Certificates and verification}\label{sec:replay}

The traces were generated with version 4.1 of Glucose~\cite{glucose} through the
\texttt{g4} backend of PySAT~\cite{pysat}\footnote{Solver interface at
\url{https://pysathq.github.io/docs/html/api/solvers.html}.}
and are stored as clausal traces~\cite{drup}. For each of the \(21\)
inputs, \texttt{drat-trim}~\cite{drattrim} verifies the derivation of
the empty clause in RUP-only mode by checking the steps it needs.
For A and R0 case~0 the checker uses \texttt{-p} to retain deleted
clauses. The other replays use the active clause database.
The \path{q7/drat-trim-replay.md} record gives the checker revision,
options, and result per input. Table~\ref{tab:checks} gives the
instance sizes.

\begin{table}[htbp]
\centering
\caption{Seven-cube verification inputs. Variables and clauses are per CNF.
Trace sizes are uncompressed totals by family.}
\label{tab:checks}
\small
\begin{tabular}{@{}lrrrr@{}}
\toprule
Statement & Inputs & Variables & Clauses & Trace size (MiB)\\
\midrule
A  & 1  & 832     & 4{,}422  & 59.0\\
R0 & 10 & 2{,}440 & 12{,}974 & 1{,}293.2\\
R1 & 10 & 3{,}438 & 14{,}978 & 2{,}512.2\\
\bottomrule
\end{tabular}
\end{table}

For a proposed derived clause \(C\), a RUP check assumes the
negation of all its literals and applies unit propagation to the
input and previously verified clauses. A conflict implies \(C\).
This test is reverse unit propagation (RUP). Inductively, every
verified clause follows from the input. Keeping earlier verified clauses
is sound by the same induction. A verified empty clause therefore proves
unsatisfiability.

The proof uses the CNF reductions above, agreement of the reconstructed
and actual inputs, and sound replay of the refutations. The encoding
audits compare clause multisets for A and ordered clause lists for R0
and R1. The supplemental PowerShell/C\# replay validated all
\(14{,}205{,}184\) additions. Its history-based checker can consult
all strictly earlier validated additions. See \path{q7/README.md}
for checker tests and recorded run times.

\section{The eight-cube}\label{sec:q8}

\subsection{Reduction to two parity profiles}

We use the following one-sided local bound. The graph on the opposite
parity side is unrestricted. The proof and its finite inputs are described
in Appendix~\ref{app:eight}.

\begin{lemma}\label{lem:q8-gate}
Let \(J\subseteq Q_6\) be \(C_4\)-free and fix a parity side \(P\).
Put \(S=\sum_{v\in P}\binom{d_J(v)}3\),
\(c=|\{v\in P:d_J(v)=6\}|\), and \(L=|\{v\in P:d_J(v)\le2\}|\).
Then \(4S+2c\le632+5L\).
\end{lemma}

\begin{proof}
Appendix~\ref{app:q8-six-face} reduces the assertion to the five
exclusions in Table~\ref{tab:eight-local}. The \texttt{q6} mode of
\path{q8/verify.py} rebuilds the \(27\) normalized inputs and
verifies their propositional refutations.
\end{proof}

Suppose that a \(C_4\)-free subgraph of \(Q_8\) has at least \(684\)
edges. Extend it by adding edges to an inclusion-maximal \(C_4\)-free
graph \(G_8\), with \(m=e(G_8)\ge684\). A missing edge at an isolated vertex
could be added without creating a square, so \(G_8\) has no isolated
vertices. On either parity side, let \(n_d\) count the vertices of
degree \(d\). Then \(\sum_d n_d=128\) and \(\sum_d d n_d=m\).

There are \(112\) six-faces. Each three-star lies in ten of them. A
vertex of degree \(d\) appears with local degree at most two in
\(\ell_d=\sum_{j=0}^2\binom dj\binom{8-d}{6-j}\) of them. Summing
Lemma~\ref{lem:q8-gate} gives
\[
  \sum_d \left(40\binom d3+2\binom d6-5\ell_d\right)n_d\le 70{,}784 .
\]
On subtracting \(402m\), adding \(1610\cdot128\), and dividing by two,
we obtain the first of the following two inequalities. The second comes
from \(\sum_d\binom d3 n_d\le1{,}792\) by Lemma~\ref{lem:star}.
\begin{equation}\label{eq:q8-budgets}
\begin{aligned}
  \sum_d C_d n_d &\le 138{,}432-201m,&
  (C_d)_{d=0}^8&=(735,534,333,177,66,0,0,105,345),\\
  \sum_d b_d n_d &\le 6{,}912-10m,&
  (b_d)_{d=0}^8&=(40,30,20,11,4,0,0,5,16).
\end{aligned}
\end{equation}
The first inequality implies \(m\le688\). To list every admissible
profile, set \(n_0=0\), enumerate \(n_1,n_2,n_3,n_4,n_7,n_8\) within
the first budget, and solve the two degree sums for \(n_5,n_6\).
Discard negative solutions and apply the second budget. For
\(m=684,\ldots,688\) this gives \(385,170,66,20,4\) profiles,
respectively. Every actual parity profile is on the list.

\subsection{Three distributions from one graph}\label{sec:q8-sampling}

Let \(\lambda\) be the distribution of \(G_8[F]\) for a uniform
four-face \(F\), and let \(\alpha\) give weight \(1/2\) to each
parity profile of \(G_8\), adding the weights if they coincide. For a
uniform ambient edge \(uv\), record whether it is retained, with bit
\(r\). For each of the seven other directions, record a state
\(s=a+2b+4c\). Here \(a,b\) are the retention bits in that direction
at \(u,v\), and \(c\) is the bit of the translated edge parallel to
\(uv\). Let \(\beta\) be the distribution of \(r\) and the
histogram of the seven sampled directions in the eight possible
states, taken up to reversing the endpoints together.
Reversal swaps states \(1,2\) and \(5,6\).

A retained edge has no state \(7\), and a missing edge has at least
one by maximality. If a retained edge has a degree-one endpoint, then
maximality forces all seven states to be \(6\) when that endpoint
is \(u\), or \(5\) when it is \(v\). Indeed, adding a missing edge
at the leaf must complete a square using its only retained edge.
For a missing edge, removing one count of state \(7\) leaves
\(\binom{13}{7}=1{,}716\) histograms of size six on eight states.
For a retained edge there are initially \(\binom{13}{6}=1{,}716\)
histograms on states \(0,\ldots,6\). Each degree-one condition
has \(\binom{10}{3}=120\) histograms, and their intersection has
\(\binom81=8\). Replacing their union by the two forced histograms
gives \(1{,}486\). Reversal fixes \(300\) of the \(3{,}202\) allowed
histograms, giving \((3{,}202+300)/2=1{,}751\) orbits.

Each distribution has mass one. Write \(x=(\lambda,\alpha,\beta)\)
and index \(\lambda\) by the entire domain of \(C_4\)-free four-cubes
up to isometry. The three distributions obey common sampling
identities. For example, if \(h_j(K)\) counts degree-\(j\) vertices in
a four-cube \(K\), choosing a vertex in a uniform four-face is
equivalent to choosing a global vertex and four of its eight directions.
Consequently
\begin{equation}\label{eq:q8-hypergeom}
 \E_\lambda\frac{h_j(K)}{16}
   =\E_\alpha\frac{\sum_d n_d\binom dj\binom{8-d}{4-j}}{128\binom84}
       \qquad(0\le j\le4).
\end{equation}
Likewise, an ambient edge together with three of its other directions
determines the same marked four-face that one obtains by choosing a
four-face and then an edge in it. Its state histogram can therefore
be counted under either \(\lambda\) or \(\beta\). The other identities
use this same sampling principle for second moments and for counts
of edge--vertex incidences. Appendix~\ref{app:eight} specifies the row
families and the local inequalities used in the certificate.

\subsection{An exact separation certificate}\label{sec:q8-separation}

\begin{lemma}[Finite separation]\label{lem:q8-separation}
There is no triple of probability distributions on the four-face,
parity-profile, and ambient-edge domains of Section~\ref{sec:q8-sampling}
that satisfies the \(356\) equalities of Table~\ref{tab:eight-rows}
and the \(83\) active inequalities of Appendix~\ref{app:eight}.
\end{lemma}

\begin{proof}
Write the identities as \(Ex=b\), with \(b_0=b_1=1\) for the
normalizations of \(\lambda,\alpha\) and all other entries zero.
Write the \(83\) inequalities used in the certificate as \(Ax\le0\).
The certificate gives integer row weights \(y,z\) with \(z\le0\)
entrywise. Set \(a=(yE+zA)/10^9\) and \(\Delta=(y_0+y_1)/10^9\).
The identities and inequalities give the lower bound
\begin{equation}\label{eq:q8-lower}
 a\cdot x\ge\Delta>10^{-4}.
\end{equation}
For an upper bound, consider each of the three distributions in \(x\)
separately. Let \(A_{\max}\) and \(B_{\max}\) be the largest coefficients
on the \(\alpha\) and \(\beta\) domains, respectively. On the
\(\lambda\) domain, the rounded four-face objective differs from the
exact coefficient by at most \(\varepsilon_4\). Its maximum \(M_4\)
is found by splitting a four-cube into two three-cubes and testing
every compatible matching between them. Each distribution has mass
one, so the exact comparisons in Table~\ref{tab:q8-separation} give
\begin{equation}\label{eq:q8-margin}
 a\cdot x\le M_4+\varepsilon_4+A_{\max}+B_{\max}<10^{-8}.
\end{equation}
This contradicts~\eqref{eq:q8-lower}. The \texttt{reconstruct} mode
of \path{q8/verify.py} forms both query files from the stored rows and
checks the uniform error bound. The \texttt{q4} mode recomputes their
maximum over all legal marked four-face completions. Appendix~\ref{app:eight}
gives the sampling, the rounding bound, and the search coverage.
\end{proof}

An inclusion-maximal counterexample yields a triple satisfying these rows.
Lemma~\ref{lem:q8-separation} therefore proves \(f(8)\le683\).
The \(682\)-edge construction in Section~\ref{sec:intro} gives the
lower bound.

\section{The asymptotic density bound}\label{sec:asymptotic}

\subsection{Two scores and the bound they imply}

We use rooted coordinate faces in the sense of Baber~\cite{baber}.
Identify a coordinate \(k\)-face \(R\) with \(\{0,1\}^k\). For a
direction \(j\) outside \(R\), its extension letter is the pair
\((m_j,o_j)\): the subset \(m_j\) of root vertices whose edges in
direction \(j\) are retained, and the retained edge set \(o_j\) on
the opposite face, transported back to \(R\) along \(j\). The root
edge set and all letters use the same labelling. An isometry of
\(R\) acts simultaneously on them. For a fixed labelled root,
independent samples of directions have independent letters.
Requiring the directions to be distinct changes
the expectation of a fixed bounded function by \(O(n^{-1})\).
The same estimate holds for any fixed number of sampled directions
and is uniform over \(G\subseteq Q_n\). All following occurrences
of the error term refer to this limit with fixed local dimensions
and coefficients.

For a \(C_4\)-free graph \(J\) on \(Q_7\), the certificate uses an
integer score \(\Gamma(J)\) and a nonnegative charge \(\mathcal A(J)\)
defined below from coordinate faces of \(J\).
Write \(S_0=2^{46}\) and \(D_0=14\cdot12\cdot10S_0\). The following
three bounds have different scopes. The first and third concern a
uniform seven-face \(F\) in an arbitrary \(C_4\)-free graph
\(G\subseteq Q_n\), while the second is pointwise on \(Q_7\).
\begin{align}
 \E\Gamma(G[F]) &\ge D_0\rho(G)-O(n^{-1}),\label{eq:asymp-score}\\
 \Gamma(J)&\le B+t\mathcal A(J),\label{eq:asymp-point}\\
 \E\mathcal A(G[F])&\le \frac{I}{32S_0}+s\rho(G)+O(n^{-1}).\label{eq:asymp-charge}
\end{align}
Here \(B\) is the score ceiling for ordinary seven-faces, \(t\ge0\)
prices the contribution of exceptional types, and \(s\) and \(I\)
control the mean charge. Table~\ref{tab:asymp-parameters} gives the
exact parameters. Combining the three estimates and letting
\(n\to\infty\) along extremal graphs gives
\begin{equation}\label{eq:asymp-result}
 \pi_\square(C_4)\le\frac{B+tI/(32S_0)}{D_0-ts}<0.59088405 .
\end{equation}
The denominator is positive since \(s<0\). The strict comparison
uses exact arithmetic with the parameters in the table.

\subsection{The seven-face score and its exceptional types}

The basic four-face scores are the edge density with rooted flag
corrections. In the limiting sampling rule, a Gram correction has
expectation equal to an average squared norm of a linear form in
the extension frequencies. The certificate also uses a quadratic
form supported on a triangle-free compatibility graph and cubic
forms given as sums of nonnegative terms. For the quadratic form,
if \(p\) is the nonnegative vector after diagonal scaling and
\(U\) is the adjacency matrix of the compatibility graph, the
Motzkin--Straus inequality gives \(p^{\mathsf T}Up\le
(\sum_i p_i)^2/2\)~\cite{ms}. Thus the quadratic form with matrix
\(\boldsymbol{1}\boldsymbol{1}^{\mathsf T}-2U\) is nonnegative.
These certificates verify the corrections as polynomials in the
extension frequencies.

To pass from four- to five-faces, the score uses paired extensions of
each three-face root. The identity
\(\|u\|^2+\|v\|^2-\|u-v\|^2=2\langle u,v\rangle\)
makes the resulting expectation a sum of squared means in the
independent sampling rule. Ceiling the unary terms and flooring the
subtracted distances only increase the score. A second paired
extension correction has the same nonnegativity property. The
integer five-face score \(\gamma_5\) satisfies
\(\E\gamma_5\ge10S_0\rho(G)-O(n^{-1})\); its formula is given in
Appendix~\ref{app:five-score}. For \(\mathcal F_k(J)\) the set of
coordinate \(k\)-faces of \(J\), define
\[
 \Gamma(J)=2\sum_{F_5\in\mathcal F_5(J)}\gamma_5(J[F_5]).
\]
The factor two counts the six-facets containing each five-face.
Uniform nested face sampling proves~\eqref{eq:asymp-score}.

\begin{lemma}[High-score catalogue]\label{lem:high-score-catalogue}
Every \(C_4\)-free \(J\subseteq Q_7\) with \(\Gamma(J)\ge B+1\)
is isometric to one of the \(192\) types in \path{frontier.json}.
Each listed type has the score recorded there.
\end{lemma}

\begin{proof}
Appendix~\ref{app:high-types} gives the anchor coverage and proves the
pruning bounds. The \texttt{all} mode of \path{asymptotic/verify.py}
starts at the four-face catalogue, regenerates the high-score five-
and six-face lists and the seven-face starting assignments, and
checks the resulting types and scores against \path{frontier.json}.
\end{proof}

\subsection{Controlling the mean charge}

For a \(k\)-face \(R\) of \(J\), let \(E_R(J)\) be its multiset of
\(q=7-k\) letters. The root term \(w_R(E_R(J))\) is zero on repeated
letters. Otherwise transport the root to its least edge mask under
cube isometries and take the lexicographically least sorted letter
list among transports with that root mask. The weight of this key
is the entry of \path{kernel.json}, with zero for absent keys.
Let \(d_4(F_4)\) be the entry of \texttt{direct\_four} in
\path{local.json} at the least edge mask of \(J[F_4]\), again with
zero for absent keys. Define
\[
 \mathcal A(J)=\sum_{k=1}^5\frac{\sum_{R\in\mathcal F_k(J)}w_R(E_R(J))}
                                  {2^{7-k}\binom7k}
       +\frac{1}{280}\sum_{F_4\in\mathcal F_4(J)}d_4(F_4).
\]

\begin{lemma}[Pointwise charge check]\label{lem:pointwise-charge}
The inequality \(\Gamma(J)\le B+t\mathcal A(J)\) holds for every
\(C_4\)-free \(J\subseteq Q_7\).
\end{lemma}

\begin{proof}
All weights are nonnegative. For \(\Gamma(J)\le B\) the inequality
follows. For the other types, the \texttt{flags} mode of
\path{asymptotic/verify.py} sums the weights over all marked roots of
each type in Lemma~\ref{lem:high-score-catalogue} and checks the inequality
in exact arithmetic.
\end{proof}

Conditional on a fixed labelled root, the expectation of its root term
is, up to \(O(n^{-1})\),
\begin{equation}\label{eq:asymp-root}
 q!\sum_{|E|=q}w_R(E)\prod_{a\in E}p_a ,
\end{equation}
where \(p_a\) is the frequency of letter \(a\). This formula keeps
all letters in the same root frame.

The certificate bounds~\eqref{eq:asymp-root} on the probability
simplex. One useful reduction is to move the combined mass of two
letters that never occur together in a positive term onto one of
them. With their sum fixed the polynomial is affine in the two
masses, so one endpoint gives at least its current value. Iterating
leaves a support in which every pair occurs together in a positive
term. Product bounds on these supports and grouping of disjoint
letter classes give the archived root inequalities. Face
averaging, including a transfer from the five-roots to their
four-facets, yields a linear combination of four-face statistics.
After adding further rooted nonnegative corrections and clearing
denominators, the certificate supplies a four-face integer price
\(h\) such that
\begin{equation}\label{eq:asymp-price}
 \E\mathcal A(G[F])-s\rho(G)\le\frac{\E h(G[F_4])}{32S_0}+O(n^{-1}),
\end{equation}
where \(F_4\) is a uniform four-face.

To bound the expectation on the right, maximize the mean of \(h\)
on the four-faces of a single six-cube. The price on one four-face
limits the prices on adjacent four-faces because their common
edges must agree. In particular, when a four-face lies outside a
distinguished list, use it as a base face of a five-cube. All eight
adjacent four-faces share a single choice for the matching in the
fifth direction. Maximizing their compatible prices with that
matching fixed improves the bound over separate maximization.
\begin{lemma}[Six-face price bound]\label{lem:six-face-price}
For every \(C_4\)-free \(L\subseteq Q_6\), the mean of \(h\) on its
coordinate four-faces is at most \(I\).
\end{lemma}

\begin{proof}
Appendix~\ref{app:charge} defines \(h\) and splits all six-faces according
to their five-facet types. The \texttt{delta} mode of
\path{asymptotic/verify.py} evaluates the four-face prices, the exact
matching bound for the five-facet cases, and the remaining six-face
enumeration. Both cases give the bound \(I\).
\end{proof}

By nested face sampling Lemma~\ref{lem:six-face-price} gives
\(\E h(G[F_4])\le I\). Substitution in~\eqref{eq:asymp-price} proves
\eqref{eq:asymp-charge} and Theorem~\ref{thm:asymptotic}.

\appendix
\raggedbottom
\section{Symmetry and encoding details}
\label{app:certificates}

\subsection{Symmetry reduction for the two missing stars}
\label{app:orbits}

Represent a \(3\)-face by a pair \((U,\xi)\), where \(U\) is its
three-coordinate free set, and \(\xi\) records its fixed bits and has
zeros on \(U\). Table~\ref{tab:orbits} uses integers whose set bits encode
\(U\) and \(\xi\).

Fix the first missing face at \((7,0)\) and write
\(J_0=\{0,1,2\}\) and \(J_1=\{3,4,5\}\). This normalization is possible by translating
from a vertex of the specified parity side of the first face and
permuting coordinates.

The stabilizer of this face preserving even
parity has independent permutations within \(J_0\) and \(J_1\) and the four
even translations supported in \(J_0\). Its order is \(3!\cdot3!\cdot4=144\).

To see that Table~\ref{tab:orbits} covers all second faces, put
\(r=|U\cap J_0|\). When \(r>0\), an even translation on \(J_0\) can
change any pattern on \(J_0\setminus U\), with one free coordinate
absorbing the parity constraint. On \(J_1\setminus U\) the invariant
is the Hamming weight of the fixed bits. Thus for \(r>0\) and fixed
weight \(w\) on \(J_1\setminus U\), the orbit has size
\[
  \binom{3}{r}^{\!2}\,2^{3-r}\binom{r}{w}.
\]
The case \(r=3,w=0\) is the first face itself. When \(r=0\),
the fixed pattern on \(J_0\) has invariant parity, giving two orbits
of size four. The listed orbits are disjoint and cover the other
\(159\) faces.

\begin{table}[htbp]
\centering
\caption{Missing-face orbits and numbers of verified proof steps.}
\label{tab:orbits}
\small
\begin{tabular}{@{}rcrrr@{}}
\toprule
Case & Second face \((U,\xi)\) & Orbit size &
\(\mathrm{R0}\) steps & \(\mathrm{R1}\) steps\\
\midrule
0 & \((7,8)\)   &  3 &    31{,}234 &   115{,}050\\
1 & \((7,24)\)  &  3 &   118{,}954 &   448{,}387\\
2 & \((7,56)\)  &  1 &   929{,}905 &   997{,}572\\
3 & \((11,0)\)  & 18 &   398{,}037 & 1{,}193{,}666\\
4 & \((11,16)\) & 36 &   556{,}394 & 1{,}096{,}140\\
5 & \((11,48)\) & 18 &   755{,}077 & 1{,}160{,}105\\
6 & \((25,0)\)  & 36 &   331{,}182 &   990{,}166\\
7 & \((25,32)\) & 36 &   600{,}382 & 1{,}164{,}830\\
8 & \((56,0)\)  &  4 &   726{,}716 & 1{,}227{,}344\\
9 & \((56,1)\)  &  4 &   391{,}579 &   680{,}982\\
\midrule
 & Total & 159 & 4{,}839{,}460 & 9{,}074{,}242\\
\bottomrule
\end{tabular}
\end{table}

Each step count is the number of verified clause additions, including
the empty clause.

\subsection{Boolean encodings}\label{app:encoding}

For an ambient edge \(e\), let \(x_e\) indicate retention in the graph.
Each coordinate square receives a clause consisting of the negations
of its four edge variables. For a candidate star with edge variables
\(x_1,x_2,x_3\) and indicator \(s\), add the three clauses
\(\neg s\vee x_i\) and the clause
\(s\vee\neg x_1\vee\neg x_2\vee\neg x_3\).
On a given parity side of a \(3\)-face there are four candidate stars.
Their pairwise exclusion clauses enforce at most one, and their
disjunction enforces at least one when required. To bound a vertex
degree between three and five, add the positive clause for each
four-element subset of its incident edges and the negative clause on
all six. For Lemma~\ref{lem:root}, use the degree bounds on the other
even vertices, six unit clauses at the origin, and exactly one even
star in each \(3\)-face. The opposite parity side is unrestricted
except for the square constraints.

For Lemma~\ref{lem:rigid}, impose the degree bounds at every vertex.
At each vertex define a degree-three indicator \(z\) using
\(\neg z\vee\bigvee_{e\in U}\neg x_e\) for each four-element subset
\(U\) of incident edges and \(z\vee\bigvee_{e\in T}x_e\) for each
three-element subset \(T\). With the degree bounds, \(z\) is true
exactly at degree three. The counter on the degree-three indicators
uses only the even vertices. For each orbit representative, forbid the
four even stars on each of the two specified missing faces and
require exactly one even star on all other faces. For each face,
the odd-side missing-star variable \(h\) is encoded by
\(h\vee s_1\vee s_2\vee s_3\vee s_4\) and the four clauses
\(\neg h\vee\neg s_i\), where the \(s_i\) are its odd star
indicators. Count exactly two true \(h\)'s and at least four
true \(z\)'s for \(\mathrm{R0}\), or exactly five true \(h\)'s
and at least eight true \(z\)'s for \(\mathrm{R1}\).

\begin{samepage}
The cardinality constraints use a sequential unary counter~\cite{sinz}.
To enforce an upper bound of \(k\) true literals among \(\ell_1,\ldots,\ell_N\), where
\(0<k<N-1\), use variables \(q_{i,j}\) for \(1\le i\le N-k\)
and \(1\le j\le k\). The four clause families are
\[
  \begin{array}{@{}rl@{\qquad}l@{}}
  \text{(i)}   & \neg\ell_i\vee q_{i,1}, &\\
  \text{(ii)}  & \neg q_{i,j}\vee q_{i+1,j}, & i<N-k,\\
  \text{(iii)} & \neg\ell_{i+j}\vee\neg q_{i,j}\vee q_{i,j+1}, & j<k,\\
  \text{(iv)}  & \neg\ell_{i+k}\vee\neg q_{i,k}. &
  \end{array}
\]
\par
\end{samepage}
Assign \(q_{i,j}\) true if the first \(i+j-1\) inputs contain
at least \(j\) true literals. This extends every assignment satisfying
the bound. Conversely, families (i)--(iii) force the
corresponding variables for all reached prefix thresholds. An induction
on the prefix length proves this claim. A threshold already reached
in the shorter prefix propagates by (ii). A threshold first reached
by the last input follows from (iii) or, for \(j=1\), from (i).
If the \((k+1)\)-st true input is at position \(p\), the forced
\(q_{p-k,k}\) and \(\ell_p\) contradict (iv). The bounds \(k=0,N-1,N\)
are represented by negative unit
clauses, one all-negative clause, and no clauses, respectively.
Lower bounds use the upper bound on the negated inputs, and exact
counts impose both bounds with separate auxiliary variables.

\section{Finite verification for the eight-cube}\label{app:eight}

\subsection{The one-sided six-face bound}\label{app:q8-six-face}

We prove Lemma~\ref{lem:q8-gate}. Lemma~\ref{lem:star} gives \(S\le160\).
Full vertices on the chosen side have pairwise distance at least four,
since two full vertices at distance two would complete a square. Five
full vertices would have total pairwise distance at least \(40\), whereas
the six coordinates contribute at most \(6\cdot6=36\). Thus \(c\le4\).
The star contributions at degrees \(3,4,5,6\) are \(1,4,10,20\), so
\(S\equiv32-L+c\pmod3\).

\begin{table}[htbp]
\caption{The five exclusions on a specified parity side of \(Q_6\).}
\label{tab:eight-local}
\centering
\small
\begin{tabular}{@{}ccclr@{}}
\toprule
 & \(L\) & \(c\) & Excluded star count & CNF inputs \\
\midrule
(i)   & 0 & any & \(S\ge159\) & 6 \\
(ii)  & 0 & 3   & \(S\ge156\) & 1 \\
(iii) & 1 & any & \(S\ge159\) & 18 \\
(iv)  & 1 & 4   & \(S\ge156\) & 1 \\
(v)   & 2 & 4   & \(S=160\)   & 1 \\
\bottomrule
\end{tabular}
\end{table}

\noindent\emph{Full vertices.}
Normalize a chosen full vertex to zero. If there are at least three
full vertices, every pair has distance four: an antipodal pair leaves
no vertex at distance at least four from both. The first two nonzero
supports then have size four and intersect in two coordinates. Write
them as \(A\cup B,A\cup C\), with \(A,B,C\) partitioning the coordinates
into pairs. A fourth full support must meet both in two coordinates and
have size four, so it is \(B\cup C\). This gives a single normalization
for either three or four full vertices.

\noindent\emph{Exclusions (i), (ii), (iv), and (v).}
In (i) the congruence leaves \((S,c)=(160,2),(159,1),(159,4)\).
The two full vertices in the first case have distance four or six.
In the case of one full vertex, the unique three-face without a star
avoids the full vertex. After normalizing that vertex to zero, this
face has one, two, or three fixed coordinates equal to one. The
remaining case uses the normalized four full vertices and allows the
hole anywhere. In (ii) and (iv) the congruence gives \(S=158\).
Exclusions (ii), (iv), and (v) each use the single normalized set of
full vertices and the indicated exact number of low vertices, including
all degrees zero, one, and two.

\noindent\emph{Exclusion (iii).}
Translate the exceptional vertex to zero. If it has degree zero,
give it one edge. If it has degree one, replace its incident edges
by a pair of directions whose two remote square edges are not both
retained. These operations preserve all other degrees on the side
and preserve \(S\). The replacement is square-free since any newly
completed square must use both edges at zero. If every pair is
blocked, all \(30\) edges between singletons and two-subsets are
retained. This case is covered by a single input with the sole edge
at zero normalized to direction zero.

In the degree-two case normalize the retained directions to \(\{0,1\}\).
If \(S=160\), a single input imposes a star in every three-face.
If \(S=159\), the congruence gives \(c=2\), and there is one hole.
Permutations preserving \(\{0,1\}\) classify it by the number \(a\) of
its free directions in the pair, the number \(b\) of its fixed ones in
the pair, and the number \(u\) of its fixed ones outside the pair.
The ranges \(0\le a\le2\), \(0\le b\le2-a\), and \(0\le u\le1+a\)
give \(16\) orbits. Permutations match each group of free, fixed-one,
and fixed-zero directions separately within the pair and its complement.

\noindent\emph{Encoding and verification.}
Use the edge bits and square clauses from Appendix~\ref{app:encoding}.
Each candidate three-star has a bit equal to the conjunction of its
three edges. A hole bit is true exactly when all candidate star bits
on the chosen side are false. Thus the hole count is \(160-S\).
Fixing a full vertex uses six unit clauses. For (iv) and (v), a low bit
at each remaining vertex is equivalent to having at most two retained
incident edges. Its two implications are cardinality bounds, and a
counter fixes the sum of low bits. There are no degree restrictions on
the opposite side. Every graph in one of the five rows yields a
satisfying assignment of the corresponding normalized input.
The \texttt{q6} mode of \path{q8/verify.py} rebuilds all \(27\) inputs
and validates their RUP refutations. The four inputs for the three
full, blocked, and two four-full cases are in
\href{https://github.com/zhelizhou/hypercube-c4/tree/\qeightrevision/q8/instances/gate}{\path{q8/instances/gate/}}.
The complete case list is in
\href{https://github.com/zhelizhou/hypercube-c4/blob/\qeightrevision/q8/proof/six-face-bound.md}{\path{q8/proof/six-face-bound.md}}.

The exclusions and the congruence give upper bounds on \(2S+c\) of
\(316,317,321,322\) when \(L=0,1,2,3\), respectively. For \(L\ge4\),
the capacity bounds give \(2S+c\le324\). Doubling each bound yields
\(4S+2c\le632+5L\), as required.

\subsection{Common rows and reconstruction}

On a marked square, one external direction gives a four-bit
connector pattern and the four edge bits of the translated
square. The five root patterns up to isometry have
\(240,212,190,187,170\) legal extension states. Both flags use
the same square labelling, and the kernel is averaged over the
stabilizer of the root. Write \(\bar Q\) for this averaged kernel
and \(s_1,s_2\) for the two states in a four-face. The local
contribution, averaged over the \(24\) marked squares, uses
\[
 \frac56\bar Q_{s_1,s_2}+
   \frac1{12}\left(\bar Q_{s_1,s_1}+\bar Q_{s_2,s_2}\right).
\]
For a square of \(G_8\) with six available directions, two draws
with replacement are distinct with probability \(5/6\).
The two diagonal terms account for a repeated direction.
The mean local score therefore exactly equals the mean of
\(v^{\mathsf T}\bar Qv/36\), where \(v\) counts the
six states at the square.

At a marked vertex the state records a rooted square up to
swapping its two directions, giving nine states. Two uniform
direction pairs in eight directions are disjoint, overlap in
one direction, or coincide with probabilities
\(15/28,12/28,1/28\), respectively. A four-face at the vertex
contains three unordered disjoint pairs, twelve overlapping
pairs, and six coincident pairs. Weighting the mean in each
group by the corresponding probability gives the exact
global second moment.

Table~\ref{tab:eight-rows} gives the common equalities in array order.
In the companion repository,
\par\nopagebreak\noindent
\href{https://github.com/zhelizhou/hypercube-c4/blob/\qeightrevision/q8/proof/certificate-rows.md}{\texttt{q8/proof/certificate-rows.md}}
gives the exact row formulas and lifts into the six blocks. We describe
the sampling behind the rows.

\begin{table}[htbp]
\caption{Equalities in the eight-cube certificate. Row indices start at zero.}
\label{tab:eight-rows}
\centering
\small
\begin{tabular}{@{}r p{0.7\linewidth}@{}}
\toprule
Rows & Equality\\
\midrule
0--1 & Normalization of the face and parity-profile distributions\\
2--6 & Local degree marginal, equation~\eqref{eq:q8-hypergeom}\\
7--70 & Symmetric second moments of edge-rooted states\\
71--88 & Degree distribution at endpoints of retained and missing edges\\
89--124 & Balances for the two root directions of a mixed vertex state\\
125--151 & Vertex state second moments grouped by local degree\\
152--355 & Histograms of three directions at an ambient edge\\
\bottomrule
\end{tabular}
\end{table}

For a histogram \(v\) with total seven, restricting to three
directions gives a count vector \(w\) with total three.
The unbiased estimator of \(vv^{\mathsf T}/49\) is
\((3ww^{\mathsf T}-2\operatorname{diag}(w))/21\).
For distinct states this follows from the inclusion probability
\(3\cdot2/(7\cdot6)\); for the same state one separates equal
directions, included with probability \(3/7\). Both endpoint
orientations are included. For a vertex of degree \(d\), the
available rooted directions are divided into \(d\) retained
and \(8-d\) missing directions. The counts of two-direction
states with \(0,1,2\) retained root edges are therefore
\(\binom{8-d}{2},d(8-d),\binom d2\). A mixed state is counted once
from each kind of direction. The endpoint row for degree \(d\)
uses \(d\,\E_\alpha n_d/512\) for retained edges and
\((8-d)\,\E_\alpha n_d/512\) for missing edges. Each is an incidence
count over \(1{,}024\) roots, with both endpoints counted.
For the grouped second moments, write \(v_i\) for the number of
pairs of state \(i\) at a vertex of \(G_8\), let \(I_g\) be the
states with \(g\) retained root edges, and put \(F_g(d)=\sum_{j\in I_g}v_j\).
The counts just stated give \(F_g(d)\) as a function of \(d\).
If \(C_{d,i}\) counts rooted occurrences of state \(i\) over the
two endpoints at an ambient edge, restricted to endpoints of
degree \(d\), then
\[
 \sum_d F_g(d)\,\E_\beta C_{d,i}
 =\tfrac12\,\E_{\text{vertex}} v_i\sum_{j\in I_g}v_j .
\]
Indeed, a given pair at a vertex is counted by both incident
root directions. The factor \(1/2\) is \(2\cdot256/1024\).
The vertex kernel sampling rule expresses the right side
by the four-face moments.
Finally, taking all \(\binom73\) choices of three directions
at a uniform edge gives the same histogram law as taking one
of the \(32\) edges of a uniform four-face. Including both
orientations gives denominators \(70\) and \(64\), respectively.
For a missing root there are \(\binom{10}{3}=120\) possible
three-direction histograms. For a retained root there are
\(\binom93=84\), giving \(204\) rows in all.

The active inequalities have the following sources. There are
\(19\) rank-one squares and \(49\) inequalities for integral
frequencies. For the latter, an integer vector \(q\) and an integer
\(u\) give the value
\(36(q\cdot v-u)(q\cdot v-u-1)\ge0\) on each count vector \(v\)
of total six. There are eight matrices on connector states,
each tested on all \(\binom{21}{6}=54{,}264\) counts of total
six in sixteen states. The minima of these eight quadratic
forms, in certificate order, are \(240,240,112,176,56,104,66,128\).
The remaining seven rows are a leaf
bound, five neighbor-degree bounds, and a conditional vertex
identity. The pure states forced at a leaf imply that its
neighbor has degree eight. If two leaves shared a neighbor,
the square along those two retained directions would force
a second incident edge at one of the leaves. For a degree-\(d\)
vertex, its neighbors have degree
sum at most \(8d-\binom d2\), since each unordered pair of its
incident retained directions forces a distinct missing
incidence at the neighbors. The active degrees are \(1,2,3,7,8\).
For the conditional identity, let \(u\) indicate the vertex
states with a retained root edge. At a degree-\(d\) vertex,
its state histogram \(v\) satisfies \(u\cdot v=28-\binom{8-d}{2}\)
and \(\boldsymbol{1}\cdot v=28\). The matrix
\(10^6(2uu^{\mathsf T}-\boldsymbol{1}\boldsymbol{1}^{\mathsf T})\)
has quadratic value
\[
 10^6\left(2\left(28-\binom{8-d}{2}\right)^2-28^2\right).
\]
The profile term cancels this quantity for each \(d=1,\ldots,8\).
Degree zero has no mass because \(G_8\) is inclusion-maximal.

The folded upper coefficients of a four-face are nonnegative.
Their total mass across the five face-root blocks is \(64{,}512\),
and their mass in the vertex block is \(64{,}512\).
The degree histogram has mass \(16\). If
\(\epsilon_{\mathrm f},\epsilon_{\mathrm v},\epsilon_{\mathrm h}\)
are the largest absolute errors on the respective coefficient
groups after rounding and rescaling, the score error is at most
\(\epsilon_{\mathrm f}+\epsilon_{\mathrm v}
 +16\epsilon_{\mathrm h}/64{,}512\). Applying this bound to the
integer weights in \texttt{EXACT\_COMBINATION.json} gives
\(\varepsilon_4\) as used in Section~\ref{sec:q8}.

\subsection{Complete maximization on four-faces}\label{app:eight-max}

Split \(Q_4\) along a fixed direction and write \(G_0,G_1\) for its
two three-cube layers in common coordinates. The joining edges
are specified by a subset \(V\subseteq\{0,1\}^3\) of their endpoints
in \(G_0\). The assembled graph is \(C_4\)-free exactly when both
layers are \(C_4\)-free and \(V\) is independent in \(G_0\cap G_1\).
Indeed a square crossing the split uses two joining edges and the
same parallel edge of both layers. This tests every square.

There are \(2{,}902\) legal labelled three-cubes and \(8{,}421{,}604\)
ordered layer pairs. Simultaneous isometries of the layers, together
with their exchange, give \(90{,}179\) pair orbits. It suffices to
try every legal \(V\) for every representative, since the same
isometry transports \(V\) and the score. The search evaluates
\(13{,}150{,}453\) such marked completions, with raw maximum
\(-37{,}440\). Dividing by \(64{,}512\cdot10^9\) gives \(M_4\) in
Table~\ref{tab:q8-separation}. The count of marked completions allows
repeated unmarked four-cube types.

\begin{table}[htbp]
\centering
\caption{Exact values for the eight-cube separation.}
\label{tab:q8-separation}
\begin{tabular}{@{}cll@{}}
\toprule
Symbol & Quantity & Exact value \\
\midrule
\(\Delta\) & Lower bound from the row combination & \(92993/500000000\) \\
\(A_{\max}\) & Maximum profile coefficient & \(276661617/1027604480000000000\) \\
\(B_{\max}\) & Maximum rooted-edge coefficient & \(1/250000000\) \\
\(\varepsilon_4\) & Uniform four-face rounding error bound & \(529874002317/601295421440000000000\) \\
\(M_4\) & Maximum rounded four-face objective & \(-13/22400000000\) \\
\bottomrule
\end{tabular}
\end{table}

The exact values give \(10^4\Delta>1\) and
\(10^8(M_4+\varepsilon_4+A_{\max}+B_{\max})<1\) by cross-multiplication.
The pinned \href{https://github.com/zhelizhou/hypercube-c4/tree/\qeightrevision/q8}{\path{q8/}}
directory contains the separation inputs. The \texttt{reconstruct} and
\texttt{q4} modes of \path{q8/verify.py} reconstruct the rational
combination and perform the complete four-face search, respectively.

\section{Exact scores and enumeration for the density bound}
\label{app:asymptotic}

The inputs are in the release asset
\href{https://github.com/zhelizhou/hypercube-c4/releases/download/q8-asymptotic-artifacts-v0.1/asymptotic-certificate.zip}{\path{asymptotic-certificate.zip}}.
Paths below start at the root of this asset. The prefix
\path{parent.zip:} identifies a file in its nested archive. The
\href{https://github.com/zhelizhou/hypercube-c4/blob/\replayrevision/asymptotic/proof.md}{\path{asymptotic/proof.md}}
note gives the mask orders, table keys, and exact price evaluators.

\begin{table}[htbp]
\centering
\caption{Exact parameters for the density bound and its finite checks.}
\label{tab:asymp-parameters}
\begin{tabular}{@{}cll@{}}
\toprule
Symbol & Role & Exact value \\
\midrule
\(B\) & Seven-face score cutoff & \(69853532380292814\) \\
\(t\) & Multiplier of the charge & \(819616459129440/819616459129441\) \\
\(s\) & Slope in the mean charge bound & \(-64043983856809\) \\
\(I\) & Upper bound on the mean four-face price & \(86291393953281818071958295771\) \\
\addlinespace
\(C_5\) & Five-face score cutoff & \(415421881252830\) \\
\(M_5\) & Upper bound on the five-face score & \(416476266812383\) \\
\(K_6\) & Six-face score cutoff & \(4989369142177746\) \\
\(M_6\) & Upper bound on the six-face score & \(4991631643675183\) \\
\addlinespace
\(U_5\) & Upper bound on the five-face mean price & \(86291393953281850225208706560\) \\
\(V_5\) & Corresponding bound outside the selected list & \(86291393953281464386203777082\) \\
\bottomrule
\end{tabular}
\end{table}

For the final comparison in~\eqref{eq:asymp-result}, the parameters
give \(10^8(B+tI/(32S_0))<59088405(D_0-ts)\) by cross-multiplication.

\subsection{Five-face score}\label{app:five-score}

Let \(f_0,f_1,f_2\) denote the three four-face scores of
Section~\ref{sec:asymptotic}. Their integer matrices are in
\path{parent.zip:q128.bin}, \path{parent.zip:q128-first.bin}, and
\path{parent.zip:q128-second.bin},
respectively, each with scale \(2^{62}\). To evaluate a score,
add \(e(F)/32\) to two mean kernel values. The first chooses a
square and a uniform labelling that sends its edge pattern to
the canonical pattern, then reads both one-direction extensions
under this labelling. The second chooses a vertex and one of the three
partitions of the four directions into pairs, then reads the two
rooted squares. The kernel entries already include the lifts
of the cubic terms. With independent direction sampling
these contributions give the Gram, copositive and cubic
polynomials described in Section~\ref{sec:asymptotic}.
In each case
\(\E f_i\ge \rho(G)-O(n^{-1})\). For a four-face \(F\), write
\(u(F)=\lceil S_0 f_0(F)\rceil+
\sum_{R\subset F,\ \dim R=3}\|b_{R,F}\|^2\). Here \(b_{R,F}\)
is the vector in \path{parent.zip:arm_factors.txt} indexed by the
isometry type of \(R\) and the least matching mask under transports
to its canonical root. This lookup uses the matching bits only.
For a five-face \(J\), set \(\theta=9604666/1073741824\).
The intermediate integer score is
\[
\begin{split}
 N(J)={}&\sum_{\substack{F\subset J\\\dim F=4}}
    \left\lceil(1-\theta)u(F)+\theta\lceil S_0 f_1(F)\rceil\right\rceil\\
 &-\sum_{\substack{R\subset J\\\dim R=3}}
    \left\lfloor(1-\theta)\|b_{R,F_i}-b_{R,F_j}\|^2\right\rfloor ,
\end{split}
\]
where \(F_i,F_j\) are the two four-faces between \(R\) and \(J\).
Each \(\|b_{R,F}\|^2\) occurs in the comparison for exactly that
three-face. On dropping the upward rounding errors, the pairwise
identity from Section~\ref{sec:asymptotic} applies. Set
\(W(J)=\sum_{\dim F=4}\lceil S_0 f_2(F)\rceil\), summing within \(J\).
For \((d,a,b)=(1048576,68033,444)\), the final score is
\[
\gamma_5(J)=
\left\lceil\frac{(d-a-b)N(J)+aW(J)+bS_0 e(J)/8}{d}\right\rceil
 +2\sum_{\substack{R\subset J\\\dim R=3}}
       \langle v_{R,F_i},v_{R,F_j}\rangle_w .
\]
The vector \(v_{R,F}\) is a second lookup, now indexed by the root
and the full extension letter, including the translated three-face.
Transport the marked extension to one of the rows of
\path{parent.zip:gamma-factors.txt} and take its \(18\)-entry vector,
with zero for an unlisted class. The inner product uses the
\(18\) nonnegative diagonal weights in that file's header.
The factors \(S_0 e(J)/8\)
and \(10S_0\rho(J)\) agree, since \(Q_5\) has \(80\) edges.
Taking expectations gives the five-face estimate used in
\eqref{eq:asymp-score}.

\subsection{Coverage above the seven-face threshold}\label{app:high-types}

Use the five-face and six-face cutoffs and upper bounds from
Table~\ref{tab:asymp-parameters}. Writing \(T=B+1\), the comparisons
that ensure the search covers every high-scoring seven-face are
\[
 2M_5+10C_5<K_6,\qquad 14K_6<T,\qquad M_6+13K_6<T .
\]
The five-face search uses the squared-norm identity to put
unary terms on the four-facets and subtract nonnegative gaps
on their intersections. Write \(u_1(F)\) for the rounded unary
summand of \(N\). The search's unary term is
\[
 \left\lceil\frac{(d-a-b)u_1(F)+a\lceil S_0 f_2(F)\rceil
                      +b S_0 e(F)/32}{d}\right\rceil
       +\sum_{R\subset F,\ \dim R=3}\|v_{R,F}\|_w^2 .
\]
At a shared three-face it subtracts the gap in \(N\) times
\((d-a-b)/d\), rounded down, and the weighted squared distance
between the two \(v\) vectors. The identity for the inner product,
\(\sum_{\dim F=4} e(F)=4e(J)\), and the rounding directions show
that this gives an upper bound on \(\gamma_5(J)\).
Let \(U_{\max}\) be the maximum unary term
over the complete four-face catalogue. Every four-facet
of a graph with \(\gamma_5>C_5\) must have unary term at least
\(C_5-9U_{\max}\), and at least one must have term at least
\(C_5/10\). Place a facet of greatest unary score as the base
and split it into two three-cubes along another direction.
The two adjacent facets sharing these layers determine
the two opposite three-cube layers and their edges to the base.
The eight remaining matching edges can be tried in every
allowed combination, checking the other four-facets and
all squares. Each true completion has one of these placements,
with any chosen isometric copy of the base. Any projected facet
bound used in pruning takes the largest unary score of all
types agreeing with the fixed bits. On a completed candidate
the score is evaluated directly.

A six-face with score at least \(K_6\) has at least three
five-facets with score above \(C_5\). At least one pair consists
of a highest facet and a nonopposite facet. A seven-face with
score at least \(T\) has at least two six-facets with scores
strictly above \(K_6\). The six-face catalogue includes equality
at \(K_6\), while the seven-face anchor condition is strict.
The two anchors can be placed either as intersecting facets,
choosing a highest facet, or as an opposite pair if the highest
facet has no other high neighbor.

At each step, mark the facets used to place an anchor in common
coordinates. For intersecting anchors the common face must have
identical edge sets under the marking. For opposite anchors
the joining edges remain to be assigned. When two side
six-facets are placed against a common bottom six-facet,
their five-dimensional intersection already has its
bottom four-face fixed. The compatibility index checks the
remaining \(80-32=48\) edge bits in that intersection.
Searching edge assignments enforces
each square by forcing its last edge absent whenever the other
three are retained. With a fixed four-facet, the program
bounds the five-face score by the maximum of \(C_5\) and the
scores of compatible catalogue entries. It compares the
joining edges and opposite four-facet in one common labelling.
When there is no fixed four-facet, the global bound \(M_5\)
is a valid cap. With at most four free edge bits the program
can instead evaluate every completion directly.

In the seven-face completion, write \(U_P\) for the cap of a
five-face \(P\) and \(c_F\) for the cap from the chosen six-facet
branch. A completion with five-face scores \(s_P\) has losses
\(\ell_P=U_P-s_P\ge0\), with \(\ell_P=0\) on completed faces,
and satisfies
\[
 \sum_{P\subset F}\ell_P\ge
 d_F=\max\left(0,\sum_{P\subset F}U_P-c_F\right).
\]
The total score is \(2\sum_P U_P-2\sum_P\ell_P\) since each
five-face lies in two six-facets. To bound the subtracted
quantity, one may sum these inequalities with nonnegative
multipliers whose sum at the two facets of each unfinished
five-face is at most two. The program uses unit multipliers
on all facets, and also uses weights two on opposite facets.
Its additional bounds put weight two on one or two facets
with no common unfinished five-face, zero on their neighbors
through unfinished faces, and one on the other facets with
positive deficit. Each of these weighted sums is valid.
Pruning occurs only when the resulting upper bound is below
the threshold. Branches that exhaust the current search budget
are carried forward with their edge assignments and bounds.

The five-face list has \(218{,}307\) marked records with score
strictly above \(C_5\). The indexed seven-face run records
\(11{,}667\) possible bottom anchors and \(2{,}592{,}200\) initial
partial assignments across the intersecting and opposite branches.
It emits \(438\) labelled outputs and closes all branches.
The replay rebuilds the start assignments by direct matching of
edges, then checks equality of the final catalogue, including
the integer score of each of the \(192\) isometry types with
score at least \(T\). The charge also uses two below-threshold
types to define the allowed patterns. Since \(\mathcal A\ge0\) and
\(T=B+1\), the \(192\) types suffice for the pointwise check
in~\eqref{eq:asymp-point}.

\subsection{Exact charge and compatibility checks}\label{app:charge}

The charge coefficients, the rooted simplex inequalities, and
the four-face price are specified by \path{transfer.json},
\path{basis.json}, \path{kernel.json}, and \path{local.json}
at the root of the certificate archive. The tables verify the root bounds for
\(390\) rooted types, \(81\) with positive terms, and transport
the five-root bounds to weighted averages over four-facets.
For example, if a group of \(q\)-element sets shares
\(q-1\) letters \(a_1,\ldots,a_{q-1}\), the sum of the
probabilities of those unordered draws is
\[
 q!\,p_{a_1}\cdots p_{a_{q-1}}\sum_{b\in B'}p_b
 \le q!/q^q ,
\]
where \(B'\) is the set of possible last letters. The inequality
is arithmetic--geometric mean on \(q\) disjoint masses of
total at most one. Paying for such a group together can give
a smaller bound than paying \(q!/q^q\) for each set.

The certificate also bounds the weights on a support of size
\(r\) by a common level plus positive remainders. The sum
of all the \(q\)-fold distinct products is at most
\(\binom rq/r^q\), by the elementary symmetric mean inequality.
For two distinct terms, put each common letter in its own
class and pair the remaining letters one from each term.
Both terms occur in the product of the \(q\) class masses, so
their combined probability is at most \(q!/q^q\).
If \(y_i\) are the individual probabilities divided by this
bound, then \(0\le y_i\le1\) and \(y_i+y_j\le1\) for distinct
terms. For nonnegative remainders \(a_i\), this implies
\(\sum_i a_i y_i\le\max(\max_i a_i,\tfrac12\sum_i a_i)\).
Indeed, if one \(y_i>1/2\), the other terms are bounded by
\(1-y_i\), and the resulting affine function is bounded by
its endpoint values on \([1/2,1]\). The groups are terms sharing
\(q-1\) letters and sets of terms that draw one letter from
each class of a partition into \(q\) disjoint classes. In the
latter case the sum of the grouped probabilities is bounded
by \(q!\) times the product of the class masses, hence again by
\(q!/q^q\). For each grouped bound the common level and group
charges are given by \texttt{group\_charges} in \path{local.json}.
The verifier applies the remainder bound after subtracting those
levels and group charges from the weights. It takes the largest
resulting conditional bound over the support cliques.

Let \(p_k(R)\) be the upper bound just computed on the
conditional root polynomial for a \(k\)-face \(R\).
The \texttt{cover} field of \path{local.json} assigns a nonnegative
value \(c_5(F)\) to each canonical four-face mask, with zero for
absent entries. On each five-root \(R\), the verifier checks that
the mean of \(c_5\) over its ten four-facets is at least \(p_5(R)\).
For a four-face \(F\), put
\[
\begin{aligned}
 U(F)={}&\frac1{32}\sum_{e\in\mathcal F_1(F)}\lceil p_1(e)\rceil
       +\frac18\sum_{K\in\mathcal F_3(F)}
          \left\lceil p_3(K)+\frac16\sum_{R\in\mathcal F_2(K)}p_2(R)\right\rceil\\
       &+\lceil p_4(F)\rceil+d_4(F)+c_5(F).
\end{aligned}
\]
Here all functions of coordinate faces are evaluated at their
induced restrictions. Uniform nested face sampling gives
\(\E\mathcal A\le\E U+O(n^{-1})\).

Let \(\lambda_0=5822002986886\).
The verifier computes \(h\) by clearing denominators in
\[
 \frac{h(F)}{32S_0}
  = U(F)+\lambda_0\frac{\lceil S_0 f_2(F)\rceil}{S_0}
    -(s+\lambda_0)\frac{e(F)}{32}+\sum_j c_j(F).
\]
Here the \(c_j\) are the further corrections in \path{local.json},
including the upward rounding used in their evaluation.
Each has nonnegative expected value up to \(O(n^{-1})\).
For example, the vertex degree moments use the nonnegative
polynomials \((8p-a)^2(8p-b)^2\) and \(p(1-p)(8p-a)^2\)
in the retained-direction probability \(p\). Their local
estimates use falling-factorial moments of the four available
directions. An edge triplet rule is checked on every histogram
of eight directions, by requiring a nonnegative mean over its
three-direction subsets. The remaining active families are
Gram and copositive root corrections and quadratic forms from
the retained and missing adjacency operators. For either edge
set, its normalized adjacency \(A_0\) and diagonal row-sum
matrix \(D_1\) give nonnegative forms \(D_1\pm A_0\), by writing
the form as a sum of squares of the values at edge endpoints
with either sign. Expanding the specified combinations gives
the corrections on small faces. The coefficients
are specified in \path{local.json}. The companion
\path{asymptotic/proof.md} lists the active correction terms and
their evaluators. The \texttt{delta} replay uses
\path{verify_full_four_exact.py} to evaluate this formula on each
four-face type. Taking expectations
gives~\eqref{eq:asymp-price}.

The complete four-face catalogue contains \(3{,}212{,}821\)
types. To certify completeness, check that every listed mask
is square-free and is the least mask in its orbit under all
\(384\) isometries. Check that these masks are strictly
increasing, and sum \(384/|\operatorname{Stab}(F)|\) over them.
The sum equals the number \(1{,}226{,}436{,}381\) of labelled
square-free four-cubes obtained by counting all compatible
joining matchings over every ordered pair of square-free
three-cube layers. Thus the price is evaluated on every
four-face type.

The \(2{,}662\) distinguished types are the zero-based indices
in \texttt{palette} in \path{local.json}, in \path{parent.zip:full_masks.bin}
order. To bound the sum of the ten prices on a five-face with a
nondistinguished facet \(F\), take
\(F\) as the base, and write \(V\subseteq\{0,1\}^4\) for the
joining matching. Each adjacent four-face has a prescribed
three-face \(R\) in \(F\) and a marked subset \(V\cap R\).
For each of the eight adjacent faces, bound the maximum
compatible price by a table indexed by \(R\) and \(V\cap R\).
Order the nondistinguished types by decreasing price, with ties
broken by increasing catalogue index. The table uses the price
at rank \(129\) in this order as a floor, so types
with lower price can be omitted when forming its entries.
Use one \(V\) for all eight, add the price of \(F\) and the
global four-face maximum for its opposite, and maximize over \(V\).
The search considers the \(128\) highest-priced nondistinguished
base types. It explicitly scans matching choices for \(62\) of them.
The rest of the \(128\) already have a base price plus nine
global maxima within the established bound. For all types
below rank \(128\), the price at rank \(129\) plus nine global
maxima suffices. This gives
a ten-face sum of at most \(10V_5\), where \(V_5\) is given in
Table~\ref{tab:asymp-parameters}.

If all the four-facets are distinguished, take a base with
maximum type index and split it into two three-cubes along
another direction.
Two neighboring four-facets determine the two opposite
three-cube layers in a consistent labelling. Enumerate the
remaining joining edges and check all squares and the other
facet types. An isometry places a maximal-index facet in any
fixed labelling of its type. The same construction one dimension
higher, using the five-face row index, covers six-faces entirely
in the selected five-face list. There are \(72{,}135\)
five-face rows. Dividing the sums by ten and rounding up gives
the bound \(U_5\). The \texttt{six\_palette} field stores the
one-based indices of the selected \(6{,}323\) five-face rows in
lexicographic order of their four-facet histograms.
The tighter bound \(V_5\) applies outside that list and includes the case of a
nondistinguished facet. A six-face with one
five-facet outside the list has mean price at most
\(\lceil(11U_5+V_5)/12\rceil\). For six-faces entirely in the
list, the \(2{,}965\) compatible rows are scanned, each summing
the twelve ten-face sums and dividing by \(120\) with upward
rounding. The exterior bound gives \(I\) in~\eqref{eq:asymp-charge},
and all enumerated rows are bounded by \(I\).
Each four-face of a six-face occurs
twice in these nested sums, so the same bound holds for the
uniform four-face mean.

The repository's
\href{https://github.com/zhelizhou/hypercube-c4/tree/\replayrevision/asymptotic}{\path{asymptotic/}}
directory describes the coefficient archive and its indexing.
The \texttt{delta} mode of \path{asymptotic/verify.py} recomputes
the refined price bound. The other modes check the rooted inequalities,
high-face enumeration, and seven-face charges.

\end{document}